\documentclass{amsart}
\usepackage{amsmath,amsfonts,amssymb,amsthm}
\usepackage{ifthen}
\usepackage{longtable}
\usepackage{array}
\usepackage{url}
\usepackage{hyperref}
\hypersetup{%
        colorlinks,
        breaklinks=true,
        plainpages=false,%
        citecolor=blue,
        linkcolor=blue,
        urlcolor=blue,
        bookmarksopen=true,%
        bookmarksnumbered=false,%
        bookmarksdepth=5%
}

\newcommand{\K}{\mathbb{K}}
\newcommand{\Z}{\mathbb{Z}}
\newcommand{\Q}{\mathbb{Q}}

\newtheorem{theorem}{Theorem}
\newtheorem{lemma}{Lemma}

\theoremstyle{remark}

\title[Narayana as product of two Fibonacci]{Narayana numbers that are products of two Fibonacci numbers}

\author[J. Odjoumani]{Japhet Odjoumani}
\address{J. Odjoumani, 
Institut de Math\'ematiques et de Sciences Physiques (IMSP),
 Universit\'e d'Abomey-Calavi (UAC), Dangbo BENIN}
\email{japhet.odjoumani@imsp-uac.org}
 
\date{\today}

\begin{document}

\begin{abstract}
Let $\{N_m\}_{m\ge0}$ be the Narayana's cows sequence given by $N_0=0$, $N_1=1=N_2=1$ and 
\[
 N_{m+3}=N_{m+2}+N_m,\quad \text{ for }\; m\geq 0
\]
and let $\{F_n\}_{n\ge0}$ be the  Fibonacci sequence. In this paper we solve explicitely the Diophantine equation 
\[
N_m=F_nF_k,
\] 
in positive unknowns $m,\,n$ and $k$. That is, we find the non-zero narayana numbers that are products of two Fibonacci numbers.
\end{abstract}
 
 
\maketitle

\section{Introduction}

The Narayana's cows sequence $\{N_m\}_{n\geq 0}$ is an integer sequence which the mathematician Narayana Pandita (see \cite{Allouche-Johnson:1996}) described as the number of cows present each year, starting from one cow in the first year, where every cow has one baby cow each year starting in its fourth year of life. It is the sequence \href{https://oeis.org/A000930}{A000930} in OEIS given by $N_0=0$, $N_1=1=N_2=1$ and 
\[
 N_{m+3}=N_{m+2}+N_m,\quad \text{ for }\; m\geq 0.
\]
Its characteristic polynomial is $X^3-X^{2}-1=(X-\alpha)(X-\beta)(X-\bar\beta)$ where 
\[
\alpha=\frac{1+r_1+r_2}{3},\;
\beta=\frac{2-(r_1+r_2)+i\sqrt{3}(r_1-r_2)}{6},
\]
$r_1=\sqrt[3]{\frac{29+3\sqrt{93}}{2}}$ and $r_2=\sqrt[3]{\frac{29-3\sqrt{93}}{2}}$.\\
For some recent studies done on Narayana's cows sequence, we refer reader to \cite{Lin:2021,Bra-Das-Guzm:2020,Bhoi-Ray:2022}. In this note we call $m$-th Narayana number the $m$-th term of Narayana's cows sequence. This is not a number $\displaystyle{N(m,\ell)=\frac{1}{m}\binom{m}{\ell}\binom{m}{\ell-1}}$ also called Narayana number (see \href{https://oeis.org/A001263}{A001263} in OEIS).\\

The Fibonacci sequence $\{F_n\}_{n\geq 0}$ is the well known sequence given by $F_0=0$, $F_1=1=F_2=1$ and 
\[
 F_{n+2}=F_{n+1}+F_{n},\quad \text{ for }\; n\geq 0.
\]
Its characteristic polynomial is $X^{2}-X-1=(X-\delta)(X-\gamma)$, where 
\[
\displaystyle{\gamma=\frac{1+\sqrt{5}}{2}\quad{\text{and}}\quad\delta=\frac{1-\sqrt{5}}{2}}.
\]
In this paper we study the diophantine equations
\begin{equation}\label{eq:main}
 N_m={F_n}{F_k},
\end{equation}
in positive unknowns $m,\,n$ and $k$. We particularly show the following. 

\begin{theorem}\label{th:main}
The only nonzero Narayana numbers that are product of two Fibonacci numbers are 
\[
 1,\,2,\,3,\,4,\,6,\;9\;\text{ and }\;13.
\]
We then see that the only nonzero Narayana numbers that are square of a Fibonacci number are $1,\;4$ and $9$.
\end{theorem}

Our method of proof involves the application of Baker's theory for linear forms in logarithms of algebraic numbers, and the Baker-Davenport reduction procedure. Computations are done with the help of a computer program in SageMath.

\section{Recalls and auxiliary Results}\label{sec:2}

\subsection{Recalls on Narayana and Fibonacci sequences}

Here we recall some properties of Narayana's cows sequence and Fibonacci sequence. Particularly the Binet formula for Narayana's cows sequence is 
\[
 N_m=a\alpha^m+b\beta^m+\bar{b}{\bar\beta}^m,\quad{\text{for integer}}\;m\ge0,
\]
where
\[
 a=\frac{\alpha}{(\alpha-\beta)(\alpha-\bar\beta)}=\frac{\alpha^2}{\alpha^3+2}\;\text{ and }\;
 b=\frac{\beta}{(\beta-\alpha)(\beta-\bar\beta)}=\frac{\beta^2}{\beta^3+2}.
\]
The minimal polynomial of $a$ over integers is $31X^3-3X-1$, with $\displaystyle{\max\left\{\left|a\right|,\,\left|b\right|\right\}<1/2}$. 

We have the numerical estimates
\begin{align*}
  &1.465<\alpha<1.466,&\nonumber\\
  &0.826<|\beta|=\left|\bar\beta\right|=\alpha^{-1/2}<0.827,&\nonumber\\
  &0.417<a<0.418,&\\
  &0.278<|b|<0.279.&\nonumber
\end{align*}
So for $m\ge1$ one proves that $e(m):=N_m-a\alpha^{m}$ satisfies
\begin{equation}\label{Naray-bound1}
 \left|e(m)\right|<0.558\alpha^{-\frac{m}{2}},
\end{equation}
and by induction 
\begin{equation}\label{Naray-bound2}
 \alpha^{m-2}\le{}N_m\le\alpha^{m-1}.
\end{equation}

The Binet formula for Fibonacci sequence is 
\[
F_n=\frac{\gamma^n-{\delta}^n}{\sqrt{5}},\quad{\text{for integer}}\;n\ge0.
\]
One has $\gamma\delta=-1$. Furthermore, for $n\ge2$ one can prove by induction that 
\begin{equation}\label{Fib-bound}
  \gamma^{n-2}\le{}F_n\le\gamma^{n-1}.
\end{equation}

We finish this subsection by noting that $\Q(\alpha)=\Q(a)$, $\Q(\alpha)\cap\Q(\gamma)=\Q$ and $3=[\Q(\alpha):\Q]\neq[\Q(\gamma):\Q]=2$. Then the numbers $\alpha, \gamma$ and $a$ are positive elements of the real field $\K=\Q(\alpha,\gamma)$ of degree $d_{\K/\Q}=6$.\\
Considering the splitting field of the polynomial $X^3-X^2-1$ over $\Q$, namely $\Q(\alpha,\beta)$, it is a Galois extension of $\Q$. In the same $\Q(\gamma,\delta)=\Q(\gamma)$, the splitting field of the polynomial $X^2-X-1$ over $\Q$, is a Galois extension of $\Q$. Then the field $\L:=\Q(\alpha,\beta,\gamma)$ is Galois extension of $\Q$ and a $\Q$-automorphism of $\L$ is for example 
\begin{equation}\label{q-automophism}
 \sigma\colon\,\alpha\mapsto\beta,\,\beta\mapsto\alpha,\,\bar\beta\mapsto\bar\beta\,\text{ and }\;x\mapsto{}x,\text{ for }\,x\in\Q(\gamma).
\end{equation}
This $\Q$-automorphism will be used later in Section \ref{sec:3}.

\subsection{Auxiliary results on linear forms in logarithms of algebraic numbers}

In this subsection, we point out some useful results from the theory of lower bounds for nonzero linear forms in logarithms of algebraic numbers. Let $\eta\neq 0$ be an algebraic number of degree $d$ and let 
\[
 a_0(X-\eta^{(1)})\cdots (X-\eta^{(d)}) \in \Z[X]
\]
be the minimal polynomial of $\eta=\eta^{(1)}$. Then the absolute logarithmic Weil height is defined by
\[
 h(\eta)=\frac{1}{d} \left(\log |a_0|+\sum_{i=1}^d \max\{0,\log|\eta^{(i)}|\}\right).
\]

This height has the following basic properties. For $\eta_1,\,\cdots,\,\eta_t$ algebraic numbers and $s\in\Z$ we have
\begin{eqnarray*}
h(\eta\pm\gamma) &\leq & h(\eta)+ h(\gamma)+\log 2,\\
h(\eta\gamma^{\pm 1}) & \leq & h(\eta) + h(\gamma),\\
h(\eta^s) &=& |s| h(\eta).
\end{eqnarray*}
In the case that $\eta$ is a rational number, say $\eta=p/q\in\Q$ with $p,\,q$ integers such that $\gcd(p,q)=1$, we have $h(p/q)=\max\{\log|p|,\log|q|\}$.\\
Now let $\K$ a real number field of degree $d_{\K}$, $\eta_1,\ldots,\eta_t \in \K$ and $b_1,\ldots,b_t\in\Z\setminus\{0\}.$ Let $B\ge\max\{|b_1|,\ldots,|b_t|\}$ and
\[
\Lambda=\eta_1^{b_1}\cdots\eta_t^{b_t}-1.
\]
Let $A_1,\ldots,A_t$ be real numbers with 
\[
A_i\ge \max\{d_{\K} h(\eta_i), |\log\eta_i|, 0.16\},\quad i=1,2,\ldots,t.
\]
With these basic notations we have the following result which is Bugeaud et al.'s version of lower bounds for linear forms in logarithms due to Matveev \cite{Matveev:2000}.
\begin{theorem}\label{Matveev}\cite[Theorem 9.4]{BMS:2006}
Assume that $\Lambda\neq 0.$ Then
\[
\log|\Lambda|>-1.4\cdot30^{t+3}\cdot t^{4.5}\cdot d_{\K}^2 \cdot(1+\log d_{\K})\cdot(1 +\log B)\cdot A_1\cdots A_t.
\]
\end{theorem}

We also need the following lemma due to Guzm\'an and Luca.
\begin{lemma}\cite[Lemma 7]{G-Luca:2014}\label{lemma_G-Luca}
If $l\ge 1,\; H>\left(4l^2\right)^l$ and $H>L/(\log L)^l,$ then 
$$ 
L<2^l H(\log H)^l.
$$
\end{lemma}

After applying these results, we find large uppers bounds for solutions of our Diophantine equation. So we use the following result of  Dujella and Peth\H{o} \cite{Dujella-Peto:1998} that is a variant of the reduction method due to Baker and Davenport \cite{Baker:1969} to reduce our bounds.

\begin{lemma}\label{Dujella-Petho:1998}
Let $M$ be a positive integer, $p/q$ be a convergent of the continued fraction expansion of the irrational number $\tau$ such that $q>6M,$ and $A, B, \mu$ be some real numbers with $A>0$ and $B>1.$ If
\[
\varepsilon:=\left\|\mu q\right\|-M\cdot\left\|\tau q\right\|>0,
\]
 then there is no solution to the inequality
\begin{equation*}
0<|u\tau-v+\mu|<AB^{-w}
\end{equation*}
in positive integers $u,v$ and $w$ with 
\[
 u\leq M\; \mbox{and}\; w\geq \frac{\log(Aq/\varepsilon)}{\log B}.
\]
\end{lemma}

\section{Proof of ours main results}\label{sec:3}

Let $(m,n,k)$ be a solution of the diophantine equation \eqref{eq:main}. We can suppose that $n\le{}k$, this is not a restriction. From \eqref{Naray-bound2} and $\eqref{Fib-bound}$, we have for $m,n\ge2$
\[
 \alpha^{m-2}<N_m=F_nF_k<\gamma^{n+k-2}\quad{\text{and}}\quad\gamma^{n+k-4}<\alpha^{m-1}.
\]
This implies that 
\begin{equation}\label{bound0-for-m}
 \frac{\log\gamma}{\log\alpha}(n+k)-2.2<m<\frac{\log\gamma}{\log\alpha}(n+k)+0.5.
\end{equation}

Furthermore, from equation \eqref{eq:main} we have 
\begin{align*}
 \left|\frac{5a\alpha^m}{\gamma^{n+k}}-1\right|&=\left|-5e(m)\gamma^{-(n+k)}-(-1)^{n}\delta^{2n}-(-1)^{k}\delta^{2k}+(-1)^{n+k}\delta^{2(n+k)}\right|\\
 &<5|e(m)|\gamma^{-(n+k)}+|\delta|^{2n}+|\delta|^{2k}+|\delta|^{2(n+k)}\\
 &<5\alpha^{\frac{-m}{2}}\gamma^{-(n+k)}+\gamma^{-2n}+\gamma^{-2k}+\gamma^{-2(n+k)}\\
 &<(5\alpha^{-1}+3)\gamma^{-2n},\quad{\text{since }}\;m\ge2,k\ge{}n.
\end{align*}
We thus obtain
\begin{equation}\label{Lamb1-upper-bound}
 \left|\frac{5a\alpha^m}{\gamma^{n+k}}-1\right|<4.91\gamma^{-2n}
\end{equation}
Putting $\displaystyle{\Lambda_1:=5a\alpha^m\gamma^{-n-k}-1}$, we have $\Lambda_1\neq 0$. Indeed, $\Lambda_1=0$ implies that $\displaystyle{5a\alpha^m=\gamma^{n+k}}$ and applying the $\Q$-automorphism $\sigma$ given in \eqref{q-automophism}, we obtain 
\begin{align*}
 &5b\beta^m=\gamma^{n+k}\\
 &5|b||\beta|^m=\gamma^{n+k}\\
 &5|b|\alpha^{-m/2}=\gamma^{n+k}.
\end{align*}
This implies that $m=-2(n+k)\frac{\log\gamma}{\log\alpha}+\log(5|b|)$ which is impossible since the first inequality in \eqref{bound0-for-m}.

We can apply now the Theorem \ref{Matveev} to $\Lambda_1$
with $t=3,$
\[
(\eta_1,\,b_1):=(\alpha,\,m),\quad (\eta_2,\,b_2):=(\gamma,\,-n-k),\quad \text{and}\quad (\eta_3,\,b_3):=(5a,1).
\]
We compute absolute logarithmic Weil height of each algebraic number and we have 
 \begin{equation}\label{eq:height}
  \begin{split}
  &h(\eta_1)=h(\alpha)<0.128, \\
  &h(\eta_2)=h(\gamma)<0.241,\\
  &h(\eta_3)=h(5a)=\log5+h(a)<2.755.
  \end{split}
 \end{equation}
We can choose 
\[
 A_1=0.768,\,A_2=1.446,\,A_3=16.53,\,B=n+k,
\]
and get
\[
 \log|\Lambda_1|>-2.742\cdot10^{14}\log(n+k).
\]
Combining with \eqref{Lamb1-upper-bound} we obtain 
\begin{align}\label{bound1-for-n}
 2n\log\gamma&<2.742\cdot10^{14}\log(n+k)+\log4.91,\nonumber\\ 
 &n\log\gamma<1.471\cdot10^{14}\log(n+k).
\end{align}

By rewriting the equation \eqref{eq:main} as 
\[
 a\alpha^m+e(m)=F_n\left(\frac{\gamma^k-\delta^k}{\sqrt{5}}\right),
\]
we have, since $n\ge2$
\begin{align*}
 \left|\frac{a\alpha^m}{F_n}-\frac{\gamma^k}{\sqrt{5}}\right|&<\frac{|e(m)|}{F_n}+\frac{|\delta|^k}{\sqrt{5}}\\
 \left|\frac{\sqrt{5}a\alpha^m}{F_n\gamma^k}-1\right|&<\frac{\sqrt{5}}{\gamma^k}\left(\frac{1}{F_n}+\frac{1}{\sqrt{5}\gamma^k}\right)\\
 &<\left(\sqrt{5}+\frac{2}{5}\right)\gamma^{-k},\quad{\text{since}}\;\gamma^k>\gamma^2>\frac{5}{2}\;{\text{and}}\;F_n\ge1,
\end{align*}
Thus putting $\Lambda_2:=\frac{\sqrt{5}a}{F_n}\alpha^m\gamma^{-k}-1$, we obtain
\begin{equation}\label{Lamb2-upper-bound}
 \left|\Lambda_2\right|<2.637\gamma^{-k}.
\end{equation}
Of course we have $\Lambda_2\neq0$.
So we can apply again the Theorem \ref{Matveev} to $\Lambda_2$ with $t=3$,
\[
 (\eta_1,\,b_1):=(\alpha,\,m),\quad (\eta_2,\,b_2):=(\gamma,\,-k),\quad \text{and}\quad (\eta_3,\,b_3):=\left(\frac{\sqrt{5}a}{F_n},1\right).
\]
We have 
\begin{align*}
 h\left(\frac{\sqrt{5}a}{F_n}\right)&\le{}h(\sqrt{5})+h(a)+\log(F_n)\le{}\frac{\log31}{3}+0.5\log(5)+\log\left(\gamma^{n-1}\right)\\
 &<2n\log\gamma.
\end{align*}
Then we can choose
\[
 A_1=0.768,\,A_2=1.446,\,A_3=12n\log\gamma,\,
 B=n+k
\]
and get
\[
 \log|\Lambda_2|>-2.018\cdot10^{14}n\log\gamma\log(n+k).
\]
Combining with \eqref{Lamb2-upper-bound} we obtain
\begin{equation}\label{bound0-for-k}
 k<4.294\cdot10^{14}n\log\gamma\log(n+k).
\end{equation}

Therefore considering the upper bound of $n\log\gamma$ from 
\eqref{bound1-for-n}, we get
\[
 k<4.294\cdot10^{14}\cdot1.471\cdot10^{14}\log^2(n+k)<6.317\cdot10^{28}\log^2(2k).
\]
Applying the Lemma \ref{lemma_G-Luca} with $l=2,\;L=2k$ and $H=1.264\cdot10^{29}$ we obtain
\begin{equation}\label{bound1-for-k}
 k<1.136\cdot10^{33}.
\end{equation}

We reduce this huge bounds by applying the Lemma \ref{Dujella-Petho:1998}. We recal that for a positive real $x$, if $\left|x-1\right|<\frac{1}{2}$ then $\left|\log{x}\right|<1.5\left|x-1\right|$ (see \cite[Lemma 4]{OZ:2021}).\\
Hence we have from \eqref{Lamb1-upper-bound},
\[
 0<\left|m\log\alpha-(n+k)\log\gamma+\log(5a)\right|<7.365\gamma^{-2n},
\]
which implies that 
\begin{equation}\label{ineq1-for-reduction}
  0<\left|m\frac{\log\alpha}{\log\gamma}-(n+k)+\frac{\log(5a)}{\log\gamma}\right|<15.306\gamma^{-2n}.
\end{equation}
Note that $\alpha$ and $\gamma$ are multiplicatively independent. Indeed, $\alpha^q=\gamma^p$ implies $2^p\alpha^q=x+y\sqrt{5}$, for some positive elements $x$ and $y$ in $\Q$. This is not possible since $3=[\Q(\alpha):\Q]\neq[\Q(\gamma):\Q]=2$ and $\gcd(2,3)=1$. Then $\displaystyle{\frac{\log\alpha}{\log\gamma}}$ is an irrational.\\
From \eqref{bound0-for-m} and \eqref{bound1-for-k} we have 
\begin{equation}\label{bound1-for-m}
 m<2.864\cdot10^{33}.
\end{equation}
So we apply Lemma \ref{Dujella-Petho:1998} with $w:=2n,$
\[
\tau:=\dfrac{\log\alpha}{\log\gamma},\quad \mu:=\dfrac{\log 5a}{\log\gamma},\quad A:=15.306,\quad B:=\gamma,\quad M:=2.864\times 10^{33}.
\]
With the help of SageMath we find that the denominator of the $72$-th convergent
\[
\dfrac{p_{72}}{q_{72}}=\dfrac{29721909555760487844132538948692737}{37417183036250693833016580755802629}
\]
of $\tau$ exceeds with $q_{72}>6M$ and $\varepsilon=0.260885028864365>0$. Thus the inequality \eqref{ineq1-for-reduction} has no solution for
\[
2n\ge \dfrac{\log (15.306\cdot{}q_{72}/\varepsilon)}{\log\gamma}\ge\dfrac{\log(15.306\cdot{}q_{72}/0.260885028864365)}{\log\gamma}\ge173.893.
\]
which implies that 
\[
n\le86.
\]
Substituting this upper bound for $n$ into \eqref{bound0-for-k}, we obtain 
\[
 k<1.778\cdot10^{16}\log(2k)
\]
We again apply Lemma \ref{lemma_G-Luca} and get
\[
 k<5.165\cdot10^{19}.
\]
From there and \eqref{bound0-for-m} we have
\[
 m<1.302\cdot10^{20}.
\]
We consider $\Lambda_2$ and we have, from \eqref{Lamb2-upper-bound}
\[
  0<\left|m\log\alpha-k\log\gamma+\log\left(\frac{\sqrt{5}a}{F_n}\right)\right|<3.956\gamma^{-k}.
\]
This implies that 
\begin{equation}\label{ineq2-for-reduction}
  0<\left|m\frac{\log\alpha}{\log\gamma}-k+\frac{\log\left(\sqrt{5}a/F_n\right)}{\log\gamma}\right|<8.22\gamma^{-k}.
\end{equation}
We then apply the Lemma \ref{Dujella-Petho:1998} with $w:=k,$
\[
\tau:=\frac{\log\alpha}{\log\gamma},\quad \mu:=\frac{\log\left(\sqrt{5}a/F_n\right)}{\log\gamma},\quad A:=8.22,\quad B:=\gamma,\quad M:=1.302\cdot10^{20}.
\]
With the help of SageMath, for $n\le86$ we find the $71$-th convergent of $\tau$
\[
 \frac{p_{71}}{q_{71}}=\frac{3194055037246978157952257926560636}{4021025019685037142147505686136939},
\]
which satisfies $q_{71}>6M$ and $\varepsilon=0.0109970619096576>0$. Hence the inequality \eqref{ineq2-for-reduction} has no solution for 
\[
 k\ge\frac{\log(8.22\cdot{}q_{71}/\varepsilon)}{\log\gamma}\ge\dfrac{\log(8.22\cdot{}q_{71}/0.0109970619096576)}{\log\gamma}\ge174.5458
\]
Thus we obtain $k\le174$ and consequently $m\le438$. We now check \eqref{eq:main} for $1\le{}n\le86,\;1\le{}k\le174$ and $1\le{}m\le438$. This is done quickly with a program on SageMath and get 
\[(m,n,k)=(m,k,n)\in\left\{
  \begin{split}
   &(1,\,1,\,1),\;(1,\,2,\,1),\;(1,\,2,\,2),\;(2,\,1,\,1),\;(2,\,1,\,2),\;(2,\,2,\,2),\\
   &(3,\,1,\,1),\;(3,\,1,\,2),\;(3,\,2,\,2),\;(4,\,1,\,3),\;(4,\,2,\,3),\;(5,\,1,\,4),\\
   &(5,\,2,\,4),\;(6,\,3,\,3),\;(7,\,3,\,4),\;(8,\,4,\,4),\;(9,\,7,\,1),\;(9,\,7,\,2).
  \end{split}\right\}
\]
This finishes the proof of the Theorem \ref{th:main}.

\section*{Acknowledgement} 
The author thanks the anonymous referee for his/her careful reading and valuable suggestions.


\begin{thebibliography}{10}

\bibitem{Allouche-Johnson:1996}
J. P. Allouche, and T. Johnson,
\newblock Narayana's cows and delayed morphisms,
\newblock {\em J. d'Inform. Music (France)}, 1996. 
\href{https://hal.archives-ouvertes.fr/hal-02986050}{hal-02986050}

\bibitem{Baker:1969}
A.~Baker and H.~Davenport,
\newblock The equations {$3x^{2}-2=y^{2}$} and {$8x^{2}-7=z^{2}$}.
\newblock {\em Quart. J. Math. Oxford Ser. (2)}, 20:129--137, 1969.

\bibitem{Bhoi-Ray:2022}
K. Bhoi, and K. P. Ray,
\newblock Fermat Numbers in Narayana's cows Sequence,
\newblock {\em Integers}, 22:\# A16, 2022.

\bibitem{Bra-Das-Guzm:2020}
J. J. Bravo, and P. Das, and S. Guzm\'an,
\newblock Repdigits in Narayana's cows sequence and their consequences,
\newblock {\em J. Integer Seq.}, 23:Article 20.8.7, 2020.

\bibitem{Bugeaud:2006}
Y.~Bugeaud, M.~Mignotte, and S.~Siksek,
\newblock Classical and modular approaches to exponential {D}iophantine
  equations. {I}. {F}ibonacci and {L}ucas perfect powers.
\newblock {\em Ann. of Math. (2)}, 163(3):969--1018, 2006.

\bibitem{Dujella-Peto:1998}
A. Dujella, and A. Peth\"o,
\newblock A generalization of a theorem of Baker and Davenport.
\newblock {\em  Quart. J. Math. Oxf. Ser. (2)}, 49:291--306, 1998.

\bibitem{G-Luca:2014}
S. Guzm\'an, and F. Luca,
\newblock Linear combinations of factorials and $s$-units in a binary recurrence sequence.
\newblock {\em  Ann. Math. Qu\'{e}}, (38):169--188, 2014.

\bibitem{Lin:2021}
X. Lin,
\newblock On the Recurrence Properties of Narayana's Cows Sequence,
\newblock {\em  Symmetry}, 149: p13, 2021. 
\href{https://doi.org/10.3390/sym13010149}{doi.org/10.3390} and \href{https://www.mdpi.com/journal/symmetry}{symmetry}

\bibitem{Matveev:2000}
E.~M. Matveev.
\newblock An explicit lower bound for a homogeneous rational linear form in
  logarithms of algebraic numbers. {II}.
\newblock {\em Izv. Ross. Akad. Nauk Ser. Mat.}, 64(6):125--180, 2000.

\bibitem{OZ:2021}
J. Odjoumani and V. Ziegler
\newblock On prime powers in linear recurrence sequences.
\newblock {\em Ann. Math. Qu\'ebec}, to appear.

\end{thebibliography}
\end{document}